\documentclass[letterpaper,10pt]{amsart}
\usepackage{amsmath,amssymb,amsxtra,amsthm}

\theoremstyle{plain}
\newtheorem{theorem}{Theorem}[section]
\newtheorem{corollary}[theorem]{Corollary}
\newtheorem{lemma}[theorem]{Lemma}
\newtheorem{proposition}[theorem]{Proposition}
\newtheorem{conjecture}[theorem]{Conjecture}

\theoremstyle{definition}
\newtheorem{definition}[theorem]{Definition}

\newtheorem{remark}[theorem]{Remark}

\newcommand{\B}{\mathbb}
\newcommand{\C}{\mathcal}

\newcommand{\ga}{\alpha}

\newcommand{\gl}{\lambda}

\newcommand{\gz}{\zeta}

\newcommand{\md}{\mbox{ \textnormal{mod} }}

\begin{document}
\title{(Non)Automaticity of Number Theoretic Functions}
\author{Michael Coons}%
\address{Simon Fraser University\\ 
Department of Mathematics\\ 
8888 University Drive\\ 
Burnaby, British Columbia V5A 1S6\\ 
Canada}

\email{mcoons@sfu.ca}
\subjclass{Primary 11J91; 11B85 Secondary 11N64}%
\keywords{automatic sequences, transcendence, Dirichlet series, multiplicative functions}%
\date{\today}

\begin{abstract} Denote by $\lambda(n)$ Liouville's function concerning the parity of the number of prime divisors of $n$. Using a theorem of Allouche, Mend\`es France, and Peyri\`ere and many classical results from the theory of the distribution of prime numbers, we prove that $\lambda(n)$ is not $k$--automatic for any $k> 2$. This yields that $\sum_{n=1}^\infty \lambda(n) X^n\in\mathbb{F}_p[[X]]$ is transcendental over $\mathbb{F}_p(X)$ for any prime $p>2$. Similar results are proven (or reproven) for many common number--theoretic functions, including $\varphi$, $\mu$, $\Omega$, $\omega$, $\rho$, and others. 
\end{abstract}

\maketitle


\section{Introduction}

In \cite{banks} it is shown that the series \begin{equation}\label{series}\sum_{n\geq 1} f(n)X^n\notin\mathbb{Z}(X)\end{equation} (is not a rational function with coefficients in $\B{Z}$) for $f$ any of the number--theoretic functions \begin{equation}\label{ntf}\varphi, \tau, \sigma, \lambda, \mu, \omega, \Omega, p, \mbox{ and } \rho.\end{equation} Here $\varphi(n)$, the Euler totient function, is the number of positive integers $m\leq n$ with $\gcd(m,n)=1$, $\tau(n)$ is the number of positive integer divisors of $n$, $\sigma(n)$ is the sum of those divisors, $\omega(n)$ is the number of distinct prime divisors of $n$, $\Omega(n)$ is the number of total prime divisors of $n$, $\lambda(n)=(-1)^{\Omega(n)}$ is Liouville's function, $\mu(n)$ is the M\"obius function defined by $$\mu(n)=\begin{cases} 1 &\mbox{if $n=1$},\\ 0 &\mbox{if $k^2|n$ for some $k\geq 2$},\\ (-1)^{\omega(n)} &\mbox{if $k^2\nmid n$ for all $k\geq 2$},\end{cases}$$ $p(n)$ is the $n$--th prime number, and $\rho(n)=2^{\omega(n)}$ counts the number of square--free positive divisors of $n$. 

In the course of this investigation we will give (or give reference to) results showing that the series $\sum_{n\geq1}f(n)X^n\in\B{Z}[[X]]$ is transcendental over $\B{Z}(X)$, for all of the functions $f$ in \eqref{ntf}. In most cases, the stronger result of transcendence of the series in $\B{F}_p[[X]]$ over $\B{F}_p(X)$ is shown. To get at these stronger results we rely upon the idea of automaticity.

Let $\mathbf{T}=(t(n))_{n\geq 1}$ be a sequence with values from a finite set. Define the {\em $k$--kernel} of $\mathbf{T}$ as the set $$\mathbf{T}^{(k)}=\{(t(k^ln+r))_{n\geq 0}:l\geq 0 \mbox{ and } 0\leq r<k^l\}.$$ Given $k\geq 2$, we say a sequence $\mathbf{T}$ is $k$--{\em automatic} if and only if the $k$--kernel of $\mathbf{T}$ is finite. Connecting automaticity to transcendence, we have the following theorem of Christol.

\begin{theorem}[Christol \cite{Chr1}] Let $\B{F}_p$ be a finite field and $(u_n)_{n\geq 0}$ a sequence with values in $\B{F}_p$. Then, the sequence $(u_n)_{n\geq 0}$ is $p$--automatic if and only if the formal power series $\sum_{n\geq 0} u_n X^n$ is algebraic over $\B{F}_p(X)$.
\end{theorem}

Since any algebraic relation in $\B{F}_p(X)$ is an algebraic relation in $\B{Z}(X)$, we have 

\begin{lemma} Let $p$ be a prime. If a series $F(X)\in\B{F}_p[[X]]$ is transcendental over $\B{F}_p(X)$ then $F(X)\in\B{Z}[[X]]$ is transcendental over $\B{Z}(X)$.
\end{lemma}

Between Allouche \cite{All1} and Yazdani \cite{Yaz1} we have that for any prime $p$, the series \eqref{series} is transcendental over $\B{F}_p(X)$ (and so over $\B{Z}(X)$ by the lemma) for $f= \varphi, \tau_k, \sigma_k,$ and $\mu$. Recall that $$\tau_k(n):=\#\{(a_1,a_2,\ldots,a_k):a_1a_2\cdots a_k=n, a_i\in\mathbb{N} \mbox{ for $i=1,\ldots,k$}\}$$ and $\sigma_k(n)$ is the sum of the $k$th powers of the divisors of $n$ (note that $\tau_2(n)=\tau(n)$ and $\sigma_1(n)=\sigma(n)$). Borwein and Coons \cite{BC1} have recently shown that the series \eqref{series} is transcendental over $\B{Z}(X)$ for any completely multiplicative function $f:\B{N}\to\{-1,1\}$ that is not identically $1$; this includes $f=\lambda$. We summarize in the following two theorems.

\begin{theorem}[Allouche \cite{All1}, Yazdani \cite{Yaz1}] The series \eqref{series} is transcendental over $\mathbb{F}_p(X)$ for $f= (g\mod v)$ with $g=\varphi, \tau_m, \sigma_m,$ and $\mu$ where $m\geq 1$ and $v\geq 2$.
\end{theorem}

\begin{theorem}[Borwein and Coons \cite{BC1}]\label{BC} The series \eqref{series} is transcendental over $\mathbb{Z}(X)$ for any nontrivial completely multiplicative function taking values in $\{-1,1\}$ (this includes $f=\lambda$).
\end{theorem}

In Section \ref{DS}, answering a question of Yazdani \cite{Yaz1}, I give the main result of this paper, the following improvement of Theorem \ref{BC}, along with many related results.

\begin{theorem}\label{nl} Liouville's function, $\lambda$, is not $k$--automatic for any $k\geq 2$, and hence $\sum_{n=1}^\infty \gl(n)X^n\in\B{F}_p[[X]]$ is transcendental over $\B{F}_p(X)$ for all $p>2$.
\end{theorem}

We can use Theorem \ref{nl} to prove the similar result for $\Omega(n)$ using the following theorem, which is a direct consequence of the definition of automaticity.

\begin{lemma}\label{Yaz} Let $t:\B{N}\to Y$ and $\Phi: Y\to Z$ be mappings. If $(t(n))_{n\geq 1}$ is $k$--automatic for some $k\geq 2$, then $(\Phi(t(n)))_{n\geq 1}$ is also $k$--automatic. 
\end{lemma}

Since $$\left(\Omega(n)\md2\right)= \frac{1-\lambda(n)}{2},$$ using Lemma \ref{Yaz} and the fact that we have $$\sum_{n\geq 1} (\Omega(n)\md 2)X^n=\sum_{n\geq 1} \Omega(n)X^n\in\B{F}_2[[X]]$$ we have the following corollary.

\begin{corollary}\label{corOmega} The function $(\Omega(n) \md 2)$ is not $2$--automatic; furthermore, the series $\sum_{n\geq 1} \Omega(n)X^n$ is transcendental over both $\B{F}_2(X)$, and $\B{Z}(X)$.
\end{corollary}

Lemma \ref{Yaz} also gives a nice corollary regarding $\tau$.

\begin{corollary} The sequence $(\tau(n) \md 2)$ is not $2$--automatic; hence the series $\sum_{n\geq 1}\tau(n)X^n$ is transcendental over both $\B{F}_2(X)$ and $\B{Z}(X)$.
\end{corollary}

\begin{proof} The function $(\tau(n)\md 2)$ is just the characteristic function of the squares, which is not $2$--automatic (see \cite{Rit1}).
\end{proof}

One of the nicest results in this area is the result of Hartmanis and Shank on the non--automaticity of the characteristic function of the primes.

\begin{theorem}[Hartmanis and Shank \cite{HS1}]\label{Min} The characteristic function of the primes, $\chi_{P}$, is not $k$--automatic for any $k\geq 2$.
\end{theorem}

In Section \ref{DS}, we give different (short and analytic--based) proofs of Theorem \ref{Min}, as well as its extension to all prime powers, and Corollary \ref{corOmega}. Many other functions are also considered in this section, such as $\rho$. As another point of interest, in Section \ref{kreg}, we address multiplicative functions which are unbounded using the generalization of $k$--automatic sequences to $k$--regular sequences.

The differences in transcendence over $\B{Z}(X)$ and $\B{F}_p(X)$ are quite pronounced. Theorem \ref{BC} gives transcendence over $\B{Z}(X)$ to a very large class of functions, many of which are $k$--automatic for some $k\geq 2$ and hence algebraic over rational functions over some finite field. For those $(f(n))_{n\geq 0}$ that are automatic, using the theory of Mahler \cite{Mahl1, Nish1} one can give transcendence results regarding the values of the series $\sum_{n=1}^\infty f(n)X^n\in\B{Z}[[X]]$. For non--automatic sequences almost no progress has been made. Indeed, it is widely believed that the number $\sum_{n=1}^\infty \gl(n)2^{-n}$ is transcendental over $\B{Q}$, and more generally we believe the following conjecture to hold, though any hope of progress is well disguised.

\begin{conjecture}\label{f2t} Let $f:\B{N}\to\{-1,1\}$ be a completely multiplicative function for which $f(p)=-1$ for at least one prime $p$. Then the number $\sum_{n=1}^\infty f(n)2^{-n}$ is transcendental over $\B{Q}$.
\end{conjecture}

\remark{As some support for this conjecture, we may focus on those sequences here which are automatic. Since all of the numbers described in Conjecture \ref{f2t} are irrational (see \cite{BC1}), by a very deep theorem of Adamczewski and Bugeaud \cite{AB1}, if for $f$ as in Conjecture \ref{f2t}, $(f(n))_{n\geq 1}$ is $k$--automatic for some $k\geq 2$, then the number $\sum_{n=1}^\infty f(n)2^{-n}$ is transcendental over $\B{Q}$.}

\section{Dirichlet Series and (non)Automaticity}\label{DS}

We rely heavily a theorem of Allouche, Mend\`es France, and Peyi\`ere \cite{All2}, and also on the details of its proof. Before proceeding to this theorem, we need some additional properties of $k$--automatic sequences (see \cite{All2} for details).

Let $k\geq 2$ and $(u(n))_{n\geq 1}$ be a $k$--automatic sequence with values in $\B{C}$. Then there exist an integer $t\geq 1$ and a sequence $(U_n)_{n\geq 1}$ with values in $\B{C}^t$ (which we denote as a column vector) as well as $k$ $t\times t$ matrices $A_1,A_2,\ldots, A_{k}$, with the property that each row of each $A_i$ has exactly one entry equal to $1$ and the rest equal to $0$ (The fact that these are $1$s and $0$s comes from the finiteness of the $k$--kernel of $(u(n))_{n\geq 1}$.), such that the first component of the vector $(U_n)_{n\geq 1}$ is the sequence $(u_n)_{n\geq 1}$ and for each $i=1,2,\ldots,k$, and for all $n\geq 1$, we have $$U_{kn+i}=A_iU_n.$$

\begin{theorem}[Allouche, Mend\`es France, and Peyi\`ere, \cite{All2}]\label{AMP} Let $k\geq 2$ be an integer and let $(u_n)_{n\geq 0}$ be a $k$--automatic sequence with values in $\B{C}$. Then the Dirichlet series $\sum_{n\geq 1}u_n n^{-s}$ is the first component of a Dirichlet vector (i.e., a vector of Dirichlet series) $G(s)$, where $G$ has an analytic continuation to a meromorphic function on the whole complex plane, whose poles (if any) are located on a finite number of left semi--lattices.
\end{theorem}

\begin{proof} We will follow the proof in \cite{All2}, but with some slight modifications. Define a Dirichlet vector $G(s)$ for $\Re s>1$ by $$G(s)=\sum_{n=1}^\infty\frac{U_n}{n^s}.$$ Since $U_{kn+j}=A_jU_n,$ we have $$G(s)=\sum_{j=1}^{k-1}\sum_{n=1}^\infty \frac{A_jU_n}{(kn+j)^s}+\sum_{n=1}^\infty \frac{A_kU_n}{(kn)^s}.$$ Writing $I$ as the $t\times t$ identity matrix, we have \begin{align*} (I-k^{-s}A_k)G(s) &= \sum_{j=1}^{k-1}\sum_{n=1}^\infty \frac{A_jU_n}{(kn+j)^s}\\
&=\sum_{j=1}^k A_j\sum_{n=1}^\infty k^{-s} n^{-s} U_n\left(1+\frac{j}{kn}\right)^{-s}\\
&=\sum_{j=1}^k A_j \sum_{m=0}^\infty {s+m-1 \choose m} (-j)^m \frac{G(s+m)}{k^{s+m}},\end{align*} and so $$(I-k^{-s}(A_0+A_1+\cdots+A_k))G(s)=\sum_{j=1}^k A_j \sum_{m=1}^\infty {s+m-1 \choose m} (-j)^m \frac{G(s+m)}{d^{s+m}}.$$ Denote $\C{A}:=k^{-1}\sum_{j=1}^k A_j$ and by $\C{M}(X)$ the transpose of the comatrix of $(\C{A}-XI)$. Multiplying the preceding equality by $\C{M}(k^{s-1})$, we have \begin{equation}\label{all3}\det(\C{A}-k^{s-1}I)G(s) = -\C{M}(k^{s-1})\sum_{j=1}^k A_j \sum_{m=1}^\infty {s+m-1 \choose m} (-j)^m \frac{G(s+m)}{k^{s+m}}.\end{equation}
For a given $s\in\B{C}$, $F(s+m)$ is bounded for $m$ large enough, so that the right--hand side of \eqref{all3} converges for $\Re s>0$ with possible poles at points $s$ for which $k^{s-1}$ is an eigenvalue of $\C{A}$. If $\Re s\in(-1,0]$, the right-hand side of \eqref{all3} converges with the possible exception of those $s$ for which $k^s$ is an eigenvalue of $\C{A}$, and so gives a meromorphic continuation of $G$ to this region with possible poles at points $s$ for which either $k^{s-1}$ or $k^s$ is an eigenvalue of $\C{A}$. Continuing this process gives an analytic continuation of $F$ to a meromorphic function on all of $\B{C}$ with possible poles at points $$s=\frac{\log \alpha}{\log k}+\frac{2\pi i}{\log k}m-l+1,$$ where $\alpha$ is an eigenvalue of $\C{A}$, $m\in\B{Z}$, $l\in\B{N}$ and $\log$ is a branch of the complex logarithm. 
\end{proof}

The beauty of this proof is in the details, which is why we have chosen to reproduce it here. Note that the possible poles are explicitly given, as is the analytic continuation. This leads to a few nice classifications regarding Dirichlet series. 

\begin{proposition} If the Dirichlet series $\sum_{n\geq 1} f(n)n^{-s}$ is not analytically continuable to the whole complex plane then $(f(n))_{n\geq 1}$ is not $k$--automatic for any $k\geq 2$.
\end{proposition}

Our first application of this is a new proof of the well--known result of Hartmanis and Shank about the non--automaticity of the characteristic function of the primes.

\begin{proof}[Proof of Theorem \ref{Min}] In 1920, Landau and Walfisz \cite{LW} proved that the Dirichlet series $P(s):=\sum_{p}p^{-s}$ is not continuable past the line $\Re s=0$. This is a consequence of the identity $$P(s)=\sum_{n\geq 1}\frac{\mu(n)}{n}\log \gz(ns).$$ Since $\gz(s)$ has a pole at $s=1$, this relationship shows that $s=1/n$ is a singular point for all square--free positive integers $n$. This sequence limits to $s=0$. Indeed, all points on the line $\Re s=0$ are limit points of the poles of $P(s)$ (see \cite[pages~215--216]{T} for details) so that this line $\Re s=0$ is a natural boundary for $P(s)$.
\end{proof}

Minsky and Papert \cite{MP1} were the first to address this question, showing that the characteristic function of the primes was not $2$--automatic. Hartmanis and Shank \cite{HS1} gave the complete result. Similar to our proof of Theorem \ref{Min}, denoting by $$\chi_{\Pi}(n):=\begin{cases} 1 & \mbox{if $n$ is a prime power}\\ 0& \mbox{otherwise},\end{cases}$$ and using the relationship $$\Pi(s):=\sum_{n\geq 1}\frac{\chi_{\Pi}(n)}{n^s}=\sum_{k\geq 1}\sum_{n\geq 1}\frac{\mu(n)}{n}\log\gz(kns),$$ we have the corresponding result for prime powers.

\begin{proposition}\label{chiPi} The sequence $(\chi_{\Pi}(n))_{n\geq 1}$ is not $k$--automatic for any $k\geq 2$.
\end{proposition}

Using Lemma \ref{Yaz} we have a result regarding $\rho$. 

\begin{proposition} Define the function $r(n)$ by $2\cdot r(n)=\rho(n)$. The sequence $(r(n) \md 2)$ is not $2$--automatic; hence $\sum_{n\geq 1}\rho(n)X^n$ is transcendental over $\B{Z}(X)$.
\end{proposition}

\begin{proof} This follows from the the fact that $(r(n)\md 2)=\chi_{\Pi}(n)$ and an application of Proposition \ref{chiPi}.
\end{proof}

As eluded to, the proof of Theorem \ref{AMP} reveals much in the way of details. Indeed, due to the explicit determination of the poles, we can can provide a very useful classification, but first, a definition. 

\begin{definition} Denote by $R(a,b;T)$ the rectangular subset of $\B{C}$ defined by $\Re s\in[a,b]$ and $\Im s\in[0,T]$, by $N_\infty(F(s),R(a,b;T))$ the number of poles of $F(s)$ in $R(a,b;T)$, and by $N_0(F(s),R(a,b;T))$ the number of zeros of $F(s)$ in $R(a,b;T)$.
\end{definition}

\begin{proposition}\label{moreT} Let $k\geq 2$, $(f(n))_{n\geq 1}$ be a $k$--automatic sequence and let $F(s)$ denote the Dirichlet series with coefficients $(f(n))_{n\geq 1}$. If $a,b\in\B{R}$ with $a<b$, then $N_\infty(F(s),R(a,b;T))=O(T)$.

Hence, if $G(s)=\sum_{n\geq 1}g(n)n^{-s}$ ($\Re s>\ga$ for some $\ga\in\B{R}$) is analytically continuable to a region containing a rectangle $R(a,b,T)$ for which $$\lim_{T\to\infty}\frac{1}{T}N_\infty(F(s),R(a,b,T))=\infty,$$ then $(g(n))_{n\geq 1}$ is not $k$--automatic for any $k\geq 2$.
\end{proposition}

\begin{proof} This is a direct consequence of the poles of $F$ being located on a finite number of left semi--lattices.
\end{proof}

From here on, we make systematic use of a classical result by von Mangoldt.

\begin{theorem}[von Mangoldt \cite{vMan1}]\label{NT} The number of zeros of $\gz(s)$ in $R(0,1;T)$ is $N_0(\gz(s),R(0,1;T))\asymp T\log T$.
\end{theorem}

\begin{theorem} The sequence $(\mu(n))_{n\geq 1}$ is not $k$--automatic for any $k\geq 2$; and hence the series $\sum_{n\geq 1} \mu(n)X^n$ is transcendental over both $\B{F}_p(X)$, for all primes $p$, and $\B{Z}(X)$.
\end{theorem}

\begin{proof} From the relationship $$\sum_{n\geq 1}\frac{\mu(n)}{n^s}=\frac{1}{\gz(s)}\quad (\Re s>1),$$ for the result, we need only show that $$\lim_{T\to\infty}\frac{1}{T}N_\infty\left(\frac{1}{\gz(s)},R(0,1:T)\right)=\infty.$$ This is given by Theorem \ref{NT}. Application of Proposition \ref{moreT} proves the theorem.
\end{proof}

It is note--worthy that our proof for $\mu(n)$ (and the proof for $|\mu(n)|$ below) does not use Cobham's theorem (see \cite{Cob2}) on rational densities: {\em if a sequence is $k$--automatic for some $k\geq 2$, then the density (provided it exists) of the occurrence of any value of that sequence is rational}. 

In a similar fashion to the above results. Using the extention to Dirichlet $L$--functions of von Mangoldt's theorem, we may generalize this result further.

\begin{lemma}\label{NTchi} We have $N_0(L(s,\chi),R(0,1;T))\asymp T\log T$.
\end{lemma}

\begin{corollary} Let $\chi$ be a Dirichlet character. Then $(\mu(n)\chi(n))_{n\geq 1}$ is not $k$--automatic for any $k\geq 2$.
\end{corollary}

\begin{proof} This follows directly from the fact that the sequence $(\mu(n)\chi(n))_{n\geq 1}$ is the sequence of coefficients of the series $\frac{1}{L(s,\chi)}$. Application of Lemma \ref{NTchi} and Proposition \ref{moreT} give the desired result.
\end{proof}

The proof of Theorem \ref{nl} rests on substantially more than the previous results of this investigation; it requires both the Prime Number Theorem (in the form below) as well as a very deep result of Selberg.

\begin{theorem}[Hadamard \cite{Had1}, de la Vall\'ee Poussin \cite{dlvp}] The Riemann zeta function has no zeros on the line $\Re s=1$.
\end{theorem}

\begin{theorem}[Selberg \cite{Sel1}] A positive proportion of the zeros of the Riemann zeta function lie on the line $\Re s=\frac{1}{2}$.
\end{theorem}

\begin{proof}[Proof of Theorem \ref{nl}] We use the identity $$\C{L}(s):=\sum_{n\geq 1}\frac{\gl(n)}{n^s}=\frac{\gz(2s)}{\gz(s)}\ \ (\Re s>1).$$ Using this identity, the poles of $\C{L}(s)$ are precisely the zeros of $\gz(s)$ that are not cancelled by the zeros of $\gz(2s)$ as well as the pole of $\gz(2s)$ at $s=\frac{1}{2}$. Selberg's theorem gives a positive proportion of zeros of $\gz(s)$ on the critical line and the Prime Number Theorem tells us that there are no zeros on the line $\Re s=1$; thus by Theorem \ref{NT}, $$N_\infty\left(\frac{\gz(2s)}{\gz(s)},R\left(\frac{1}{2},\frac{1}{2};T\right)\right)\asymp T\log T.$$ Application of Proposition \ref{moreT} gives the result.
\end{proof}

Invoking a stronger form of Selberg's theorem, we may include many more number--theoretic functions in our investigation.

\begin{theorem}[Conrey \cite{Con1}] More than two--fifths of the zeros of the Riemann zeta function lie on the critical line.
\end{theorem}

Conrey's theorem gives the following corollary.

\begin{corollary}\label{CC} Less than three--tenths of the zeros of the Riemann zeta function lie on any line $\Re s\neq \frac{1}{2}$.
\end{corollary}

\begin{proof} Recall that if $\gz(s)=0$, then by the functional equation $\gz(1-s)=0$. The corollary then follows from the elementary observation that $2\cdot\frac{3}{10}+\frac{2}{5}=1$.
\end{proof}

\begin{theorem}\label{kfree} For $k\geq 2$, the functions $q_m(n)$ $(m\geq 2)$ are not $k$--automatic; and hence $\sum_{n\geq 1}q_m(n)X^n$ for each $m\geq 2$ is transcendental over both $\B{F}_p(X)$, for all primes $p$, and $\B{Z}(X)$.
\end{theorem} 

\begin{proof} Note the identities for $\Re s>1$: $$\sum_{n\geq 1}\frac{q_m(n)}{n^s}=\frac{\gz(s)}{\gz(ms)}\qquad (m\geq 2).$$ Our result relies on $\gz(s)/\gz(mz)$ (for each $m\geq 2$) having more than $O(T)$ poles in some rectangle. Corollary \ref{CC} gives $$N_\infty\left(\frac{\gz(s)}{\gz(ms)},R\left(\frac{1}{2m},\frac{1}{2m};T\right)\right)\asymp T\log T.$$ Application Proposition \ref{moreT} finish the proof.
\end{proof}

Recall that for $m\geq 2$ $$q_m(n)=\begin{cases}0 & \mbox{if $p^m|n$ for any prime $p$}\\ 1 & \mbox{otherwise;}\end{cases}$$ hence $|\mu(n)|=q_2(n)$, so that the function $|\mu(n)|$ is provided for in the above corollary.

\section{Dirichlet Series and (non)Regularity}\label{kreg}

Taking the definition from \cite{AS}, we say that a sequence $\mathbf{S}:=(s(n))_{n\geq 0}$ taking values in a $\B{Z}$--module $R$ is a {\em $k$--regular sequence} (or just $k$--regular) provided there exist a finite number of sequences over $R$, $\{(s_1(n))_{n\geq 0},\ldots,(s_s(n))_{n\geq 0}\}$, with the property that every sequence in the $k$--kernel of $\mathbf{S}$ is a $\B{Z}$--linear combination of the $s_i$; that is, $\mathbf{S}$ is $k$--regular provided the $k$--kernel of $\mathbf{S}$ is finitely generated (as opposed to being finite in the case of $k$--automatic).

Using this definition, let $k\geq 2$ and $(s(n))_{n\geq 1}$ be a $k$--regular sequence with values in $\B{C}$. Then similar to the automatic case, there exist an integer $t\geq 1$ and a sequence $(V_n)_{n\geq 1}$ with values in $\B{C}^t$ (which we denote as a column vector) as well as $k$ $t\times t$ matrices $B_1,B_2,\ldots, B_{k}$ with integer entries (no longer just $1$s and $0$s as in the automatic case), such that the first component of the vector $(V_n)_{n\geq 1}$ is the sequence $(v_n)_{n\geq 1}$ and for each $i=1,2,\ldots,k$, and for all $n\geq 1$, we have $$V_{kn+i}=B_iV_n.$$ These properties give the analogue of Theorem \ref{AMP} to $k$--regular sequences.

\begin{theorem} Let $k\geq 2$ be an integer and let $(v_n)_{n\geq 0}$ be a $k$--regular sequence with values in $\B{C}$. Then the Dirichlet series $\sum_{n\geq 1}v_n n^{-s}$ is the first component of a Dirichlet vector (i.e., a vector of Dirichlet series) $G(s)$, where $G$ has an analytic continuation to a meromorphic function on the whole complex plane, whose poles (if any) are located on a finite number of left semi--lattices.
\end{theorem}

The proof of this theorem is exactly that of Theorem \ref{AMP} with $V_i$ and $B_i$ substituted for $U_i$ and $A_i$, respectively, for each $i$.

We now have the same useful corollaries that we had for $k$--automatic sequences.

\begin{corollary}\label{moreTreg} Let $k\geq 2$. The following properties hold:
\begin{itemize}
\item[(i)] If the Dirichlet series $\sum_{n\geq 1} f(n)n^{-s}$ is not analytically continuable to the whole complex plane then $(f(n))_{n\geq 1}$ is not $k$--regular.
\item[(ii)] If $G(s)=\sum_{n\geq 1}g(n)n^{-s}$ ($\Re s>\ga$ for some $\ga\in\B{R}$) is analytically continuable to a region containing a rectangle $R(a,b,T)$ for which $$\lim_{T\to\infty}\frac{1}{T}N_\infty(G(s),R(a,b,T))=\infty,$$ then $(g(n))_{n\geq 1}$ is not $k$--regular.
\end{itemize}
\end{corollary}

\begin{theorem} The function $\varphi(n)$ is not $k$--regular for any $k\geq 2$.
\end{theorem}

\begin{proof} From the relationship $$\sum_{n\geq 1}\frac{\varphi(n)}{n^s}=\frac{\gz(s-1)}{\gz(s)}\qquad (\Re s>2),$$ and the lack of zeros of $\gz(s-1)$ in the region $0\leq\Re s\leq 1$ as given by the Prime Number Theorem, we need only show that $$\lim_{T\to\infty}\frac{1}{T}N_\infty\left(\frac{1}{\gz(s)},R(0,1:T)\right)=\infty.$$ This is given by Theorem \ref{NT}. Application of the Corollary \ref{moreTreg} proves the theorem.
\end{proof}

\begin{theorem} For $k\geq 2$, the functions $\rho(n)$, $\tau(n^2)$ and $\tau^2(n)$ are not $k$--regular.
\end{theorem} 

\begin{proof} Note the identities for $\Re s>1$: $$\sum_{n\geq 1}\frac{\rho(n)}{n^s}=\frac{\gz^2(s)}{\gz(2s)},\qquad \sum_{n\geq 1}\frac{\tau(n^2)}{n^s}=\frac{\gz^3(s)}{\gz(2s)},\qquad \sum_{n\geq 1}\frac{\tau^2(n)}{n^s}=\frac{\gz^4(s)}{\gz(2s)}.$$ Since a multiple zero from the numerator is only counted once, our result relies on $\gz(s)/\gz(2z)$ having more than $O(T)$ poles in some rectangle. This follows directly from the proof of Theorem \ref{kfree}.
\end{proof}

\begin{theorem} The functions $\omega(n)$ and $\Omega(n)$ are not $k$--regular for any $k\geq 2$.
\end{theorem}

\begin{proof} This follows from the proof of Theorem \ref{Min} and the identities $$\sum_{n\geq 1} \frac{\omega(n)}{n^s}=\gz(s)\sum_{k\geq 1} \frac{\mu(k)}{k}\log \gz(ks),\  \mbox{and}\  \sum_{n\geq 1} \frac{\Omega(n)}{n^s}=\gz(s)\sum_{k\geq 1} \frac{\varphi(k)}{k}\log \gz(ks)$$ and the added stipulation that there are no zeros of $\gz(s)$ on the line $\Re s=0$; this is provided for by the Prime Number Theorem and the symmetry of zeros of the Riemann zeta function about the critical line as given by the functional equation for $\gz(s)$.
\end{proof}

Some of these results can be found from another direction using our knowledge of their non--automaticity and the following theorem (see Chapter 16 of \cite{AS} for details).

\begin{theorem}[Allouche and Shallit \cite{AS}] If the integer sequence $(f(n))_{n\geq 0}$ is $k$--regular, then for all integers $m\geq 1$, the sequence $(f(n)\md m)_{n\geq 0}$ is $k$--automatic.
\end{theorem}

Thus if there exists an $m\geq 1$ for which $(f(n)\md m)_{n\geq 0}$ is not $k$--automatic, then $(f(n))_{n\geq 0}$ is not $k$--regular. Hence the results of the previous sections give non--regularity results for each of $\omega$, $\Omega$, $\tau$, and $\rho$. It is also worth noting that a sequence is $k$-regular and takes on only finitely many values if and only if it is $k$--automatic (again, see \cite{AS}). This provides a nice relationship for non--regularity results for characteristic functions like $q_m$ ($m\geq 2$), $\chi_P$, and $\chi_{\Pi}$.

\section{Concluding Remarks}

There is much to do in this area, and it seems that the available methods and results leave many ideas ripe for development.

Concerning transcendence of power series for these functions, one need not dig so deeply to give transcendence results over $\B{Z}(X)$ or $\B{Q}(X)$ using theorems like the following.

\begin{theorem}[Fatou \cite{Fat1}]\label{fatou} If $F(X)=\sum_{n\geq 1}f(n)X^n\in\B{Z}[[X]]$ converges inside the unit disk, then either $F(X)\in\B{Q}(X)$ or  $F(X)$ is transcendental over $\B{Q}(X)$.
\end{theorem}

Carlson \cite{Car1}, proving a conjecture of P\'olya, added to Fatou's theorem.

\begin{theorem}[Carlson \cite{Car1}] A series $F(X)=\sum_{n\geq 1}f(n)X^n\in\B{Z}[[X]]$ is either rational or it admits the unit circle as a natural boundary.
\end{theorem}

Recall that if $f(n)=O(n^d)$ for some $d$, the series $F(X)=\sum_{n\geq 1}f(n)X^n$ $\in\B{Z}[[X]]$ has the unit circle as a natural boundary, so that by the combination of the above two theorems of Carlson and Fatou, such a series is transcendental over $\B{Q}(X)$. This gives very quick transcendence results for series $F(X)$ with $f(n)=\varphi(n), \tau(n^2), \tau^2(n), \omega(n),$ and $\Omega(n)$. Noting that by the Prime Number Theorem, $p(n)\sim n\log n=O(n^2)$, we have the following result for the $n$th prime number.

\begin{proposition} The series $\sum_{n\geq 1}p(n)X^n\in\B{Z}[[X]]$ is transcendental over $\B{Q}(X)$, and hence also over $\B{Z}(X)$.
\end{proposition}

The ideas of $k$-regularity may be exploitable to give transcendence results using the following theorem of Allouche and Shallit from \cite{AS} and a combination of the above theorems in this section, though it seems at this point that a case by case analysis would be necessary, which we believe would not make for easy reading. 

\begin{theorem}[Allouche and Shallit \cite{AS}] Let $K$ be an algebraically closed field (e.g., $\B{C}$). Let $(s(n))_{n\geq 0}$ be a sequence with values in $K$. Let $S(X)=\sum_{n\geq 0}s(n)X^n$ be a formal power series in $K[[X]]$. Assume that $S$ represents a rational function of $X$. Then $(s(n))_{n\geq 0}$ is $k$--regular if and only if the poles of $S$ are roots of unity.
\end{theorem}

One may be able to form this into more rigid and inclusive theorems and as such, this seems a worthy endeavor.

Concerning more specific functions, the non-automaticity of $\gl(n)$ (and similarly $\mu(n)$) is somewhat weak compared to the expected properties of the correlation. One expects that for any $A,B,a,b\in\B{N}$ with $aB\neq Ab$ $$\left|\sum_{n\leq x}\gl(An+B)\gl(an+b)\right|=o(x),$$ so that not only should the $k$--kernel be infinite (as shown in this paper), but no two sequences of $\gl$ on distinct arithmetic progressions should be equal. In this sense, the Liouville function should be a sort of ``worst case scenario'' for non--automaticity concerning multiplicative functions. It would be worthwhile to develop this idea further. Some results in this vein are known. Indeed, Hartmanis and Shank \cite{HS1} have shown that the primes can be recognized by a linearly growing automaton, but they cannot be given by logarithmically growing one.

\remark{All zeta quotient identities as well as the properties of the Riemann zeta function that were used in this paper can be found in Titchmarsh's monograph \cite{T}.}




\providecommand{\bysame}{\leavevmode\hbox to3em{\hrulefill}\thinspace}
\providecommand{\MR}{\relax\ifhmode\unskip\space\fi MR }
\providecommand{\MRhref}[2]{%
  \href{http://www.ams.org/mathscinet-getitem?mr=#1}{#2}
}
\providecommand{\href}[2]{#2}

\end{document}